\documentclass[3p,11pt]{elsarticle}
\usepackage{amssymb,amsfonts,amsmath,natbib}
\usepackage{lmodern}

\graphicspath{{figures/en/}}

\journal{Journal Title}

\numberwithin{equation}{section}

\begin{document}

\begin{frontmatter}

%% title, authors, affilations
%\input{chapters/en/Title}
\title{On exact solutions of some important nonlinear conformable time-fractional differential equations}
\author{Erdo\u{g}an Mehmet \"Ozkan$^{\rm a}$
\vspace{6pt} \\\vspace{6pt}  $^{a}${\em{Yildiz Technical University, Department of Mathematics, Istanbul, Turkey;}} \\ Corresponding author: mozkan@yildiz.edu.tr.
}

%% abstract
\begin{abstract}
The nonlinear fractional Boussinesq equations are known as the fractional differential equation class that has an important place in mathematical physics. In this study, a method called $\Big(\frac{G'}{G^2}\Big)$-extension method which works well and reveals exact solutions is used to examine nonlinear Boussinesq equations with conformable time-fractional derivatives. This method is a very useful approach and extremely useful compared to other analytical methods. With the proposed method, there are three unique types of solutions such as hyperbolic, trigonometric and rational solutions. This approach can similarly be applied to other nonlinear fractional models.
\end{abstract}

%% keywords
\begin{keyword}
Conformable fractional derivative,  $\Big(\frac{G'}{G^2}\Big)$-expansion method \sep Exact solutions 
\MSC[2010] 26A33 \sep 35R11 \sep 83C15
%% or \MSC[2008] code \sep code (2000 is the default)
\end{keyword}

\end{frontmatter}
\section{Introduction}

Nonlinear partial differential equations (NPDEs) are very useful equations that are important in different fields including engineering sciences, mathematics, fluid dynamics, and physics. Many real world problems have been modeled thanks to NPDEs. Numerous different and robust mathematical methods have been used to obtain exact solutions of NPDEs \cite{Inan, Yil1, Yil2, Yil3, Abdou, Wakil, Wazwaz,Zayed, Manafian, Kudryashov, Akbulut}. Fractional calculus is also a new field that has grown in popularity over the past few decades. Various physical phenomena are made easily solvable in fractional differential equations (FDEs), such as viscoelasticity, plasma, solid mechanics, optical fibers, signal processing, electromagnetic waves, fluid dynamics, biomedical sciences and diffusion processes. Researchers have made these equations attractive by using various methods to obtain exact solutions of FDEs. Some of these methods and the important ones in the literature can be listed as $\Big(\frac{G'}{G}\Big)$-expansion method \cite{Bekir3}, the improved extended tanh-coth method \cite{Gomez}, Lie group analysis method \cite{Chen1}, first integral method \cite{Lu}, exp-function method \cite{Zheng}, sub-equation method \citep{Aksoy}, functional variable method \cite{Matinfar}.

Boussinesq type equations can be considered as the first model for nonlinear, dispersive wave propagation and describe the surface water waves whose horizontal scale is much larger than the depth of the water \cite{Madsen}. They can be noted as a critical class of fractional differential equations in mathematical physics. Recently, different techniques have been used to find analytical and numerical solutions of Boussinesq equations. These can be mentioned as invariant subspace method \cite{Sahadevan}, a newly developed method called the expansion method \cite{Hosseini1}, exp-function method \cite{Yaslan}. 

Conformable time fractional Boussinesq equations, which can be used to describe the propagation of the wave in the magnetic field, describes surface water waves. In addition, it is an important nonlinear model that occurs in physics, hydromechanics and optics. It is also known that it can be used to describe the physical direction of wave propagation in plasma and nonlinear wave \cite{Triki, Abazari, Abdelgawad, Biswas5, Biswas6}. The coupled Boussinesq equations occur in shallow water waves for bilayer fluid flow. This happens when there is an accidental oil spill from a ship, which causes a slick of oil to float above the water slide \cite{Jawad}.

The study includes five chapters. In the first chapter, fractional calculations are mentioned along with brief information about nonlinear partial and fractional differential equations. In addition, Boussinesq equations and coupled Bousinesq equations that form the basis of the publication are explained. In the second chapter, definition and theories of comformable fractional derivative are given. In the third chapter $\Big(\frac{G'}{G^2}\Big)$-expansion method which is an important method for the solution  is explained. In Chapter four, the application of finding exact solutions of comformable time fractional Boussinesq equation and coupled conformable time fractional Boussinesq equations with this method has been made. Finally, an explanation of the results is given in the last section.

\section{Conformable fractional derivative and its properties}
\noindent
Various definitions have been introduced for the fractional derivative such as Riemann-Liouville, Gr\"unwald-Letnikov and Caputo. All fractional derivatives provide linearity, but these definitions have some flaws. Khalil et al. \cite{Khalil} described a new, simple fractional derivative named the comformable fractional derivative that overcomes these shortcomings. The conformable derivative has the following definitions and properties as given in \cite{Khalil, Abdel, Eslami}.

\textbf{Definition}: Suppose $f:[0,\infty)\rightarrow \mathbf{R}$ be a function. The $\alpha$th-order conformable fractional derivative of $f$ is defined by
$$ _t D^\alpha(f)(t)=\lim_{\varepsilon \to 0} \frac{f(t+\varepsilon t^{1-\alpha})-f(t)}{\varepsilon},\ \ \  t>0, \ \ \ \alpha \in (0,1). $$

\textbf{Theorem1}: Suppose $\alpha \in (0, 1]$ and $f, g$ be $\alpha$-differentiable at $t>0$. Then, 
\begin{description}
\item \ \ \ (1) $_t D^\alpha(af+bg)=aT_\alpha(f)+bT_\alpha(g), \ \ \ \forall a,b \in \mathbf{R}.$
\item \ \ \ (2) $_t D^\alpha(t^{p})=pt^{p-\alpha}, \ \ \ \forall p \in \mathbf{R}.$
\item \ \ \ (3) $_t D^\alpha(\lambda)=0$ ($\lambda$ is any constant).
\item \ \ \ (4) $_t D^\alpha(fg)=f(_t D^\alpha(g))+g(_t D^\alpha(f)).$
\item \ \ \ (5) $_t D^\alpha(\frac{f}{g})=\frac{g(_t D^\alpha(f))-f(_t D^\alpha(g))}{g^{2}}.$
\item \ \ \ (6) If $f$ is differentiable, then $_t D^\alpha(f)(t)=t^{1-\alpha}\frac{df}{dt}.$
\end{description}

\textbf{Theorem2}: Let $f, g:(0, \infty)\rightarrow  \mathbf{R}$ be $\alpha$-differentiable functions, when $0<\alpha \le 1$ . Then, 
$$ _t D^\alpha(f \circ g)(t)=t^{1-\alpha}g'(t)f'(g(t)). $$

\section{Algorithm of $\Big(\frac{G'}{G^2}\Big)$-expansion method}
\noindent
In this section, we give an presentation of $\Big(\frac{G'}{G^2}\Big)$-expansion method for solving a general space time fractional differential equation as 
\begin{equation} \label{eq1}
H\Big(u, \frac{\partial^{\alpha}u}{\partial t^\alpha}, \frac{\partial u}{\partial x}, \frac{\partial^{2\alpha}u}{\partial t^{2\alpha}},\frac{\partial^2 u}{\partial x^2}...\Big)=0,
\end{equation}
where $u$ is an unknown function and $H$ is a polynomial of $u$ and its partial fractional derivatives \cite{Chen}. It can be introduced as follows.
\vspace{0.3cm}

\textbf{Step1}: By applying the transformation
\begin{equation} \label{eq2}
\begin{aligned}
& u(x,t)=U(\varepsilon), \\
& \varepsilon=x+k\frac{t^{\alpha}}{\alpha},
\end{aligned}
\end{equation}
which $k$ are nonzero constant coefficients, Eq.\eqref{eq1} can be reduced  into the ordinary differential equation
\begin{equation}\label{eq3}
Q\Big(U,\frac{dU}{d\varepsilon}, \frac{d^{2}U}{d\varepsilon^{2}}, \frac{d^{3}U}{d\varepsilon^{3}}, ...\Big)=0.
\end{equation}
\vspace{0.3cm}

\textbf{Step2}: Assume that Eq.\eqref{eq3} has the solution
\begin{equation}\label{eq4}
U(\varepsilon)=a_{0}+\sum_{i=1}^{M}\Bigg[{a_{i}\Bigg(\frac{G'}{G^2}\Bigg)^i+b_{i}\Bigg(\frac{G'}{G^2}\Bigg)^{-i}}\Bigg],
\end{equation}
and
\begin{equation}\label{eq5}
\Bigg(\frac{G'}{G^2}\Bigg)' =\mu+\lambda\Bigg(\frac{G'}{G^2}\Bigg)^2
\end{equation}
where $\lambda \neq 0$ $\mu \neq 1$ are integers and $a_{0}$, $a_{i}, b_{i} (i=1,2,...,M)$ are real constants to be determined. The balancing number $M$ is a positive integer which can be determined balancing the highest derivative terms with the highest power nonlinear terms in Eq.\eqref{eq3}. 
\vspace{0.3cm}

\textbf{Step3}: Substituting \eqref{eq4} into \eqref{eq3} using \eqref{eq5} and setting the coefficients of all powers of $\frac{G'}{G^2}$ to zero, the system of algebraic equations will be obtained. This algebraic system is solved by using programs such as Maple to obtain the values of the unknown constants $a_{0}$, $a_{i}$ and $b_{i} (i=1,2,...,M)$. More clearly, the value of M can be found from Eq. \eqref{eq3} as follows where $D(U(\varepsilon))=M$ is degree of $U(\varepsilon)$.
\begin{equation*} 
\begin{aligned}
& D\bigg[\frac{d^{q}U}{d\varepsilon^{q}}\bigg]=M+q, \\
& D\bigg[U^{r}\bigg(\frac{d^{q}U}{d\varepsilon^{q}}\bigg)^{s}\bigg]=Mr+s(q+M).
\end{aligned}
\end{equation*}
\vspace{0.3cm}

\textbf{Step4}: On the basis of the general solution of Eq. \eqref{eq5}, the required exact solutions will be obtained in the following three cases: 
\vspace{0.3cm}

\noindent
\emph{Case 1.} If $\lambda\mu>0$, then
\begin{equation} \label{eq6}
\frac{G'}{G^2}=\sqrt{\frac{\mu}{\lambda}}\Bigg(\frac{C\cos(\sqrt{\mu\lambda}\varepsilon)+D\sin(\sqrt{\mu\lambda}\varepsilon)}{D\cos(\sqrt{\mu\lambda}\varepsilon)-C\sin(\sqrt{\mu\lambda}\varepsilon)}\Bigg),
\end{equation}
\vspace{0.3cm}

\noindent
\emph{Case 2.} If $\lambda\mu<0$, then
\begin{equation} \label{eq7}
\frac{G'}{G^2}=-\frac{\sqrt{|\mu\lambda|}}{\lambda}\Bigg(\frac{C\sinh(2\sqrt{|\mu\lambda|}\varepsilon)+C\cosh(2\sqrt{|\mu\lambda|}\varepsilon)+D}{C\sinh(2\sqrt{|\mu\lambda|}\varepsilon)+C\cosh(2\sqrt{|\mu\lambda|}\varepsilon)-D}\Bigg),
\end{equation}
\vspace{0.3cm}

\noindent
\emph{Case 3.} If $\lambda \neq 0$ and $\mu=0$, then
\begin{equation} \label{eq8}
\frac{G'}{G^2}=-\frac{C}{\lambda(C\varepsilon+D)},
\end{equation}
where $C,D$ are nonzero constants.

\section{Applications}

\noindent
In this section, two important nonlinear fractional differential equations will be solved by $\Big(\frac{G'}{G^2}\Big)$-expansion method.
\subsection{Conformable time fractional Boussinesq equation}
\noindent
Let us consider the following nonlinear equation \cite{Hosseini}:
\begin{equation} \label{eq9}
\frac{\partial^{2\alpha}u}{\partial t^{2\alpha}}-\frac{\partial^2 u}{\partial x^2}-\frac{\partial^2 u^2}{\partial x^2}+\frac{\partial^4 u}{\partial x^4}=0, \ \ \ \   0<\alpha\leq1.
\end{equation}

\noindent
Applying the transformation Eq.\eqref{eq2} to this equation, it is obtained the following ordinary differential equation
\begin{equation} \label{eq10}
(k^2-1)\frac{d^2 U}{d\varepsilon^2}-\frac{d^2 (U^2)}{d\varepsilon^2}+\frac{d^4 U}{d\varepsilon^4}=0.
\end{equation}
Integrating Eq.\eqref{eq10} twice and taking constants of integration as zero, the following equation is obtained.
\begin{equation} \label{eq11}
\frac{d^2 U}{d\varepsilon^2}+(k^2-1)U-U^2=0.
\end{equation}
Balancing the terms of $U^2$ and $\frac{d^2 U}{d\varepsilon^2}$ in Eq. \eqref{eq11}, we get $M=2$. For $M=2$, solution will get the form
\begin{equation}\label{eq12}
U(\varepsilon)=a_{0}+a_{1}\Big(\frac{G'}{G^2}\Big)+a_{2}\Big(\frac{G'}{G^2}\Big)^2+b_{1}\Big(\frac{G'}{G^2}\Big)^{-1}+b_{2}\Big(\frac{G'}{G^2}\Big)^{-2},
\end{equation}
where $a_{0}, a_{1}, a_{2}, b_{1}$ and $b_{2}$ are unknown constants.

\noindent
Substituting \eqref{eq12} using \eqref{eq5} into \eqref{eq11}, setting the coefficients of all powers of $\Big(\frac{G'}{G^2}\Big)$ to zero, nonlinear algebraic equations are achieved. The occuring algebraic system is solved by aid of Maple to find the values of unknown constants.
\begin{equation*}
\begin{aligned}
& \Big(\frac{G'}{G^2}\Big)^{-4} &:& \ \ \ \ 6\mu^2b_{2}-b_{2}^2=0, \\
& \Big(\frac{G'}{G^2}\Big)^{-3} &:& \ \ \ \ 2\mu^2b_{1}-2b_{1}b_{2}=0, \\
& \Big(\frac{G'}{G^2}\Big)^{-2} &:& \ \ \ \  8\mu\lambda b_{2}+(k^2-1)b_{2}-b_{1}^2-2b_{2}a_{0}=0, \\
& \Big(\frac{G'}{G^2}\Big)^{-1} &:& \ \ \ \ 2\mu\lambda b_{1}+(k^2-1)b_{1}-2b_{2}a_{1}-2b_{1}a_{0}=0, \\
& \Big(\frac{G'}{G^2}\Big)^{0} &:& \ \ \ \ 2\lambda^2 b_{2}+2\mu^2 a_{2}
+(k^2-1)a_{0}-a_{0}^2-2b_{2}a_{2}-2b_{1}a_{1}=0, \\
& \Big(\frac{G'}{G^2}\Big)^{1} &:& \ \ \ \  2\mu\lambda a_{1}+(k^2-1)a_{1}-2b_{1}a_{2}+2a_{1}a_{0}=0, \\
& \Big(\frac{G'}{G^2}\Big)^{2} &:& \ \ \ \ 8\mu\lambda a_{2}+(k^2-1)a_{2}-a_{1}^2-2a_{2}a_{0}=0, \\
& \Big(\frac{G'}{G^2}\Big)^{3} &:& \ \ \ \ 2\lambda^2a_{1}-2a_{1}a_{2}=0, \\
& \Big(\frac{G'}{G^2}\Big)^{-4} &:& \ \ \ \ 6\lambda^2a_{2}-a_{2}^2=0.
\end{aligned}
\end{equation*}

\noindent
Solving this system of equations through Maple, we get the following results;
\vspace{0.3cm}

\noindent
\textbf{Set1}: $k=\mp \sqrt{16\lambda\mu+1}, \ \ a_{0}=12\lambda\mu, \ \ a_{1}=b_{1}=0, \ \ a_{2}=6\lambda^2, \ \ b_{2}=6\mu^2$,

\noindent
\textbf{Set2}: $k=\mp \sqrt{4\lambda\mu+1}, \ \ a_{0}=6\lambda\mu, \ \ a_{1}=b_{1}=b_{2}=0, \ \ a_{2}=6\lambda^2$,

\noindent
\textbf{Set3}: $k=\mp \sqrt{4\lambda\mu+1}, \ \ a_{0}=6\lambda\mu, \ \ a_{1}=b_{1}=a_{2}=0, \ \ b_{2}=6\mu^2$,

\noindent
\textbf{Set4}: $k=\mp \sqrt{1-4\lambda\mu}, \ \ a_{0}=2\lambda\mu, \ \ a_{1}=b_{1}=b_{2}=0, \ \ a_{2}=6\lambda^2$,

\noindent
\textbf{Set5}: $k=\mp \sqrt{1-4\lambda\mu}, \ \ a_{0}=2\lambda\mu, \ \ a_{1}=b_{1}=a_{2}=0, \ \ b_{2}=6\mu^2$,

\noindent
\textbf{Set6}: $k=\mp \sqrt{1-16\lambda\mu}, \ \ a_{0}=-4\lambda\mu, \ \ a_{1}=b_{1}=0, \ \ a_{2}=6\lambda^2, \ \ b_{2}=6\mu^2$.

\vspace{0.3cm}
\noindent
The above set of values gives the following exact solutions for conformable time fractional Boussinesq equation

\vspace{0.3cm}
\textbf{Solution 1}: 

\vspace{0.2cm}
\noindent
(i) If $\lambda\mu>0$, the trigonometric solution is found 
\begin{equation*}
\begin{aligned}
& u_{1}(x,t)=U_{11}(\varepsilon)=12\lambda\mu+6\lambda^2\Bigg(\sqrt{\frac{\mu}{\lambda}}\Bigg(\frac{C\cos(\sqrt{\mu\lambda}\varepsilon)+D\sin(\sqrt{\mu\lambda}\varepsilon)}{D\cos(\sqrt{\mu\lambda}\varepsilon)-C\sin(\sqrt{\mu\lambda}\varepsilon)}\Bigg)\Bigg)^2 \\
& \ \ \ \ \ \ \ \ \ \ \ \ \ \ \ \ \ \ \ \ \ \ \ \ \ \ \ \ \ \ \ \ \ +6\mu^2\Bigg(\sqrt{\frac{\mu}{\lambda}}\Bigg(\frac{C\cos(\sqrt{\mu\lambda}\varepsilon)+D\sin(\sqrt{\mu\lambda}\varepsilon)}{D\cos(\sqrt{\mu\lambda}\varepsilon)-C\sin(\sqrt{\mu\lambda}\varepsilon)}\Bigg)\Bigg)^{-2}.
\end{aligned}
\end{equation*}

\noindent
(ii) If $\lambda\mu<0$, the hyperbolic solution is obtained
\begin{equation*}
\begin{aligned}
& u_{2}(x,t)=U_{12}(\varepsilon)=12\lambda\mu+6\lambda^2\Bigg(-\frac{\sqrt{|\mu\lambda|}}{\lambda}\Bigg(\frac{C\sinh(2\sqrt{|\mu\lambda|}\varepsilon)+C\cosh(2\sqrt{|\mu\lambda|}\varepsilon)+D}{C\sinh(2\sqrt{|\mu\lambda|}\varepsilon)+C\cosh(2\sqrt{|\mu\lambda|}\varepsilon)-D}\Bigg)\Bigg)^2 \\
& \ \ \ \ \ \ \ \ \ \ \ \ \ \ \ \ \ \ \ \ \ \ \ \ \ \ \ \ \ \ \ \ \ +6\mu^2\Bigg(-\frac{\sqrt{|\mu\lambda|}}{\lambda}\Bigg(\frac{C\sinh(2\sqrt{|\mu\lambda|}\varepsilon)+C\cosh(2\sqrt{|\mu\lambda|}\varepsilon)+D}{C\sinh(2\sqrt{|\mu\lambda|}\varepsilon)+C\cosh(2\sqrt{|\mu\lambda|}\varepsilon)-D}\Bigg)\Bigg)^{-2}.
\end{aligned}
\end{equation*}

\noindent
(iii) If $\lambda\neq0,\ \ \ \mu=0$, the rational solution is found
\begin{equation*}
\begin{aligned}
& u_{3}(x,t)=U_{13}(\varepsilon)=6\Bigg(\frac{C}{C\varepsilon+D}\Bigg)^2
\end{aligned}
\end{equation*}

\noindent
In (i)-(iii), $\varepsilon=x \mp \sqrt{16\lambda\mu+1}\big(\frac{t^{\alpha}}{\alpha}\big) $.

\vspace{0.3cm}
\textbf{Solution 2}: 

\vspace{0.2cm}
\noindent
(i) If $\lambda\mu>0$, the trigonometric solution is found 
\begin{equation*}
u_{4}(x,t)=U_{21}(\varepsilon)=6\lambda\mu+6\lambda^2\Bigg(\sqrt{\frac{\mu}{\lambda}}\Bigg(\frac{C\cos(\sqrt{\mu\lambda}\varepsilon)+D\sin(\sqrt{\mu\lambda}\varepsilon)}{D\cos(\sqrt{\mu\lambda}\varepsilon)-C\sin(\sqrt{\mu\lambda}\varepsilon)}\Bigg)\Bigg)^2.
\end{equation*}

\noindent
(ii) If $\lambda\mu<0$, the hyperbolic solution is obtained
\begin{equation*}
u_{5}(x,t)=U_{22}(\varepsilon)=6\lambda\mu+6\lambda^2\Bigg(-\frac{\sqrt{|\mu\lambda|}}{\lambda}\Bigg(\frac{C\sinh(2\sqrt{|\mu\lambda|}\varepsilon)+C\cosh(2\sqrt{|\mu\lambda|}\varepsilon)+D}{C\sinh(2\sqrt{|\mu\lambda|}\varepsilon)+C\cosh(2\sqrt{|\mu\lambda|}\varepsilon)-D}\Bigg)\Bigg)^{2}.
\end{equation*}

\noindent
(iii) If $\lambda\neq0,\ \ \ \mu=0$, the rational solution is found
\begin{equation*}
\begin{aligned}
& u_{6}(x,t)=U_{33}(\varepsilon)=6\Bigg(\frac{C}{C\varepsilon+D}\Bigg)^2
\end{aligned}
\end{equation*}

\noindent
In (i)-(iii), $\varepsilon=x \mp \sqrt{4\lambda\mu+1}\big(\frac{t^{\alpha}}{\alpha}\big) $.

\vspace{0.3cm}
\textbf{Solution 3}: 

\vspace{0.2cm}
\noindent
(i) If $\lambda\mu>0$, the trigonometric solution is found 
\begin{equation*}
u_{7}(x,t)=U_{31}(\varepsilon)=6\lambda\mu+6\mu^2\Bigg(\sqrt{\frac{\mu}{\lambda}}\Bigg(\frac{C\cos(\sqrt{\mu\lambda}\varepsilon)+D\sin(\sqrt{\mu\lambda}\varepsilon)}{D\cos(\sqrt{\mu\lambda}\varepsilon)-C\sin(\sqrt{\mu\lambda}\varepsilon)}\Bigg)\Bigg)^{-2}.
\end{equation*}

\noindent
(ii) If $\lambda\mu<0$, the hyperbolic solution is obtained
\begin{equation*}
u_{8}(x,t)=U_{32}(\varepsilon)=6\lambda\mu+6\mu^2\Bigg(-\frac{\sqrt{|\mu\lambda|}}{\lambda}\Bigg(\frac{C\sinh(2\sqrt{|\mu\lambda|}\varepsilon)+C\cosh(2\sqrt{|\mu\lambda|}\varepsilon)+D}{C\sinh(2\sqrt{|\mu\lambda|}\varepsilon)+C\cosh(2\sqrt{|\mu\lambda|}\varepsilon)-D}\Bigg)\Bigg)^{-2}.
\end{equation*}

\noindent
In (i)-(ii), $\varepsilon=x \mp \sqrt{4\lambda\mu+1}\big(\frac{t^{\alpha}}{\alpha}\big) $.

\vspace{0.3cm}
\textbf{Solution 4}:

\vspace{0.2cm}
\noindent
(i) If $\lambda\mu>0$, the trigonometric solution is found 
\begin{equation*}
u_{9}(x,t)=U_{41}(\varepsilon)=2\lambda\mu+6\lambda^2\Bigg(\sqrt{\frac{\mu}{\lambda}}\Bigg(\frac{C\cos(\sqrt{\mu\lambda}\varepsilon)+D\sin(\sqrt{\mu\lambda}\varepsilon)}{D\cos(\sqrt{\mu\lambda}\varepsilon)-C\sin(\sqrt{\mu\lambda}\varepsilon)}\Bigg)\Bigg)^2 
\end{equation*}

\noindent
(ii) If $\lambda\mu<0$, the hyperbolic solution is obtained
\begin{equation*}
u_{10}(x,t)=U_{42}(\varepsilon)=2\lambda\mu+6\lambda^2\Bigg(-\frac{\sqrt{|\mu\lambda|}}{\lambda}\Bigg(\frac{C\sinh(2\sqrt{|\mu\lambda|}\varepsilon)+C\cosh(2\sqrt{|\mu\lambda|}\varepsilon)+D}{C\sinh(2\sqrt{|\mu\lambda|}\varepsilon)+C\cosh(2\sqrt{|\mu\lambda|}\varepsilon)-D}\Bigg)\Bigg)^2 
\end{equation*}

\noindent
(iii) If $\lambda\neq0,\ \ \ \mu=0$, the rational solution is found
\begin{equation*}
\begin{aligned}
& u_{11}(x,t)=U_{43}(\varepsilon)=6\Bigg(\frac{C}{C\varepsilon+D}\Bigg)^2
\end{aligned}
\end{equation*}

\noindent
In (i)-(iii), $\varepsilon=x \mp \sqrt{1-4\lambda\mu}\big(\frac{t^{\alpha}}{\alpha}\big) $.

\vspace{0.3cm}
\textbf{Solution 5}: 

\vspace{0.2cm}
\noindent
(i) If $\lambda\mu>0$, the trigonometric solution is found 
\begin{equation*}
u_{12}(x,t)=U_{51}(\varepsilon)=2\lambda\mu+6\mu^2\Bigg(\sqrt{\frac{\mu}{\lambda}}\Bigg(\frac{C\cos(\sqrt{\mu\lambda}\varepsilon)+D\sin(\sqrt{\mu\lambda}\varepsilon)}{D\cos(\sqrt{\mu\lambda}\varepsilon)-C\sin(\sqrt{\mu\lambda}\varepsilon)}\Bigg)\Bigg)^{-2}.
\end{equation*}

\noindent
(ii) If $\lambda\mu<0$, the hyperbolic solution is obtained
\begin{equation*}
u_{13}(x,t)=U_{52}(\varepsilon)=2\lambda\mu+6\mu^2\Bigg(-\frac{\sqrt{|\mu\lambda|}}{\lambda}\Bigg(\frac{C\sinh(2\sqrt{|\mu\lambda|}\varepsilon)+C\cosh(2\sqrt{|\mu\lambda|}\varepsilon)+D}{C\sinh(2\sqrt{|\mu\lambda|}\varepsilon)+C\cosh(2\sqrt{|\mu\lambda|}\varepsilon)-D}\Bigg)\Bigg)^{-2}.
\end{equation*}

\noindent
In (i)-(ii), $\varepsilon=x \mp \sqrt{1-4\lambda\mu}\big(\frac{t^{\alpha}}{\alpha}\big) $.

\vspace{0.3cm}
\textbf{Solution 6}: 

\vspace{0.2cm}
\noindent
(i) If $\lambda\mu>0$, the trigonometric solution is found 
\begin{equation*}
\begin{aligned}
& u_{14}(x,t)=U_{61}(\varepsilon)=-4\lambda\mu+6\lambda^2\Bigg(\sqrt{\frac{\mu}{\lambda}}\Bigg(\frac{C\cos(\sqrt{\mu\lambda}\varepsilon)+D\sin(\sqrt{\mu\lambda}\varepsilon)}{D\cos(\sqrt{\mu\lambda}\varepsilon)-C\sin(\sqrt{\mu\lambda}\varepsilon)}\Bigg)\Bigg)^2 \\
& \ \ \ \ \ \ \ \ \ \ \ \ \ \ \ \ \ \ \ \ \ \ \ \ \ \ \ \ \ \ \ \ \ +6\mu^2\Bigg(\sqrt{\frac{\mu}{\lambda}}\Bigg(\frac{C\cos(\sqrt{\mu\lambda}\varepsilon)+D\sin(\sqrt{\mu\lambda}\varepsilon)}{D\cos(\sqrt{\mu\lambda}\varepsilon)-C\sin(\sqrt{\mu\lambda}\varepsilon)}\Bigg)\Bigg)^{-2}.
\end{aligned}
\end{equation*}

\noindent
(ii) If $\lambda\mu<0$, the hyperbolic solution is obtained
\begin{equation*}
\begin{aligned}
& u_{15}(x,t)=U_{62}(\varepsilon)=-4\lambda\mu+6\lambda^2\Bigg(-\frac{\sqrt{|\mu\lambda|}}{\lambda}\Bigg(\frac{C\sinh(2\sqrt{|\mu\lambda|}\varepsilon)+C\cosh(2\sqrt{|\mu\lambda|}\varepsilon)+D}{C\sinh(2\sqrt{|\mu\lambda|}\varepsilon)+C\cosh(2\sqrt{|\mu\lambda|}\varepsilon)-D}\Bigg)\Bigg)^2 \\
& \ \ \ \ \ \ \ \ \ \ \ \ \ \ \ \ \ \ \ \ \ \ \ \ \ \ \ \ \ \ \ \ \ +6\mu^2\Bigg(-\frac{\sqrt{|\mu\lambda|}}{\lambda}\Bigg(-\frac{C\sinh(2\sqrt{|\mu\lambda|}\varepsilon)+C\cosh(2\sqrt{|\mu\lambda|}\varepsilon)+D}{C\sinh(2\sqrt{|\mu\lambda|}\varepsilon)+C\cosh(2\sqrt{|\mu\lambda|}\varepsilon)-D}\Bigg)\Bigg)^{-2}.
\end{aligned}
\end{equation*}

\noindent
(iii) If $\lambda\neq0,\ \ \ \mu=0$, the rational solution is found
\begin{equation*}
\begin{aligned}
& u_{16}(x,t)=U_{63}(\varepsilon)=6\Bigg(\frac{C}{C\varepsilon+D}\Bigg)^2
\end{aligned}
\end{equation*}

\noindent
In (i)-(iii), $\varepsilon=x \mp \sqrt{1-16\lambda\mu}\big(\frac{t^{\alpha}}{\alpha}\big) $.
\vspace{0.2cm}

\noindent
In Figure 1, the physical properties of Eq. \eqref{eq9} whose solutions are used in \textbf{Solution 2} have been shown for $\alpha = 0.5$ with some special values.
\begin{figure}
\centering
  \includegraphics[width=14cm]{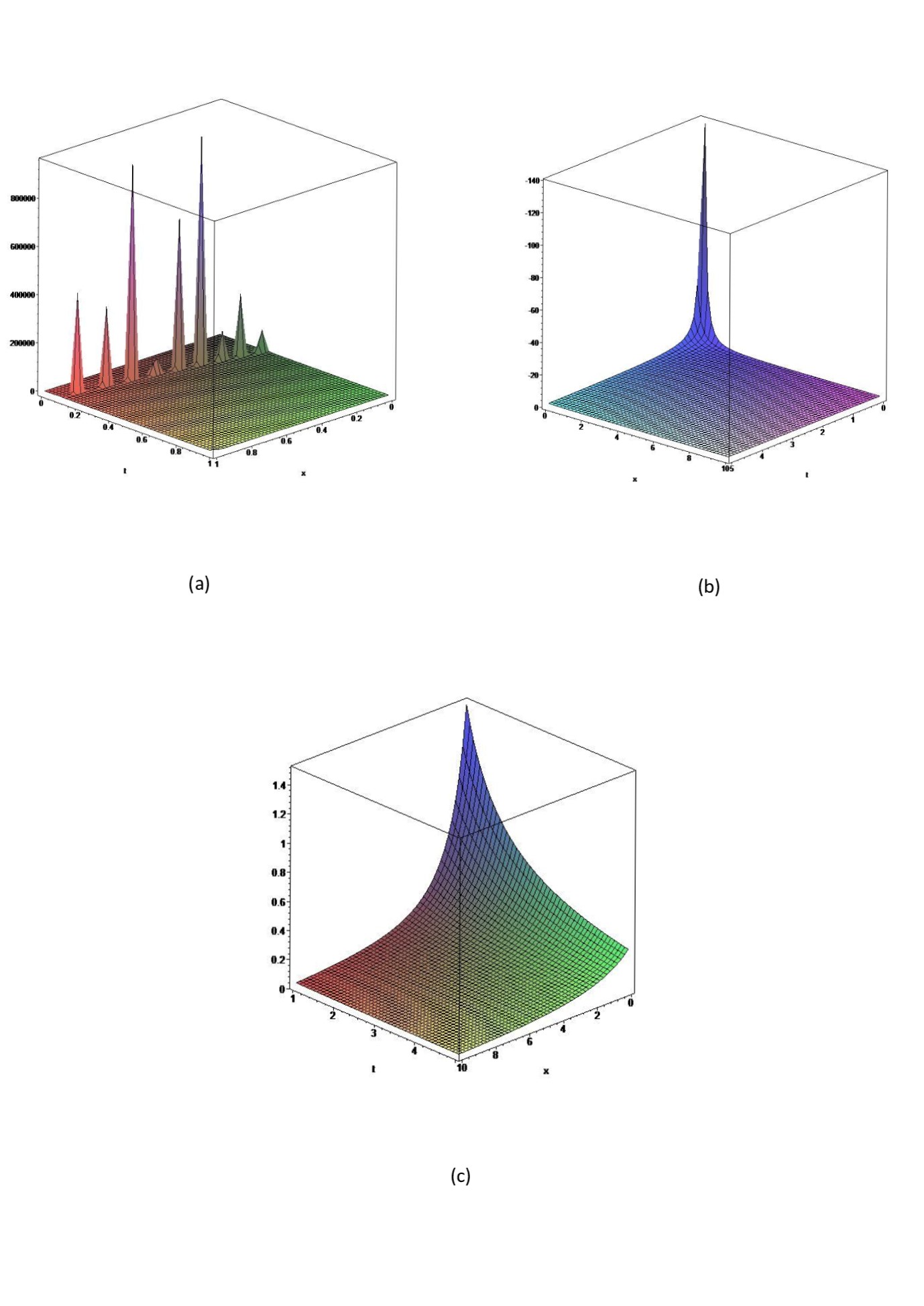}\\
  \caption{(a) The trigonometric solution of $u(x,t)$ when $ \lambda=\mu=C=D=1.$ (b) The hyperbolic solution of $u(x,t)$ when $ \lambda=0.5, \mu=-0.3, C=D=1$. (c) The rational solution of $u(x,t)$ when $ \lambda=1, \mu=0, C=D=1$.}
\end{figure}

\subsection{Coupled conformable time fractional  Boussinesq equations}

\noindent
Let us consider the following nonlinear equation system\cite{Hosseini}:
\begin{equation}
 \begin{aligned}
   & \frac{\partial^{\alpha}u}{\partial t^{\alpha}}+\frac{\partial v}{\partial x}=0,  \label{eq:1a} \\ 
   & \frac{\partial^{\alpha}v}{\partial t^{\alpha}}+\beta\frac{\partial}{\partial x}(u^2)-\gamma\frac{\partial^3 u}{\partial x^3}=0,
\end{aligned}
\end{equation}
where the time fractional derivative is $\alpha$th order and $0 <\alpha\leq 1$.

\noindent
Applying the transformation  Eq.\eqref{eq2} to this system, it is found the following ordinary differential equation system
\begin{equation}\label{eq13}
 \begin{aligned}
  & k\frac{dU}{d\varepsilon}+\frac{dV}{d\varepsilon}=0, \\
  & k\frac{dV}{d\varepsilon}+\beta\frac{d}{d\varepsilon}(U^2)-\gamma \frac{d^3U}{d\varepsilon^3}=0.
\end{aligned}
\end{equation}

\noindent
Integrating Eq.\eqref{eq13} once and taking constant of integration as zero, the following equations are found:

\begin{subequations}
    \begin{align}
     & kU+V=0,  \label{eq:2a} \\ 
     & kV+\beta U^2- \gamma \frac{d^2U}{d\varepsilon^2} =0  \label{eq:2b}.
\end{align}
\end{subequations}

\noindent
From Eq.\eqref{eq:2a}, we get
\begin{equation}\label{eq14}
V=-kU.
\end{equation}

\noindent
Substituting Eq.\eqref{eq14} into Eq. \eqref{eq:2b}, it is obtained the following differential equation
\begin{equation}\label{eq15}
-k^2 U+\beta U^2-\gamma \frac{d^2U}{d\varepsilon^2}=0.
\end{equation}

\noindent
Balancing the terms of $U^2$ and $\frac{d^2 U}{d\varepsilon^2}$ in Eq. \eqref{eq15}, we get $M=2$. For $M=2$, solution will obtain the form

\begin{equation}\label{eq16}
U(\varepsilon)=a_{0}+a_{1}\Big(\frac{G'}{G^2}\Big)+a_{2}\Big(\frac{G'}{G^2}\Big)^2+b_{1}\Big(\frac{G'}{G^2}\Big)^{-1}+b_{2}\Big(\frac{G'}{G^2}\Big)^{-2},
\end{equation}
where $a_{0}, a_{1}, a_{2}, b_{1}$ and $b_{2}$ are unknown constants.

\noindent
Substituting \eqref{eq16} using \eqref{eq5} into \eqref{eq15}, setting the coefficients of all powers of $\Big(\frac{G'}{G^2}\Big)$ to zero, nonlinear algebraic equations are obtained. The occurring algebraic system is solved by help of Maple to find the values of unknown constants.
\begin{equation*}
\begin{aligned}
& \Big(\frac{G'}{G^2}\Big)^{-4} &:& \ \ \ \ \beta b_{2}^2-6\gamma \mu^2b_{2}=0, \\
& \Big(\frac{G'}{G^2}\Big)^{-3}&:& \ \ \ \ 2 \beta b_{1}b_{2}-2 \gamma \mu^2b_{1}=0, \\
& \Big(\frac{G'}{G^2}\Big)^{-2}&:& \ -k^2b_{2}+\beta b_{1}^2+2 \beta b_{2}a_{0}-8\gamma \mu \lambda b_{2}=0, \\
& \Big(\frac{G'}{G^2}\Big)^{-1}&:& \ -k^2b_{1}+2\beta b_{2}a_{1}+2\beta b_{1}a_{0}-2\gamma \mu \lambda b_{1}=0, \\
& \Big(\frac{G'}{G^2}\Big)^{0}&:& \ -k^2a_{0}+\beta a_{0}^2+2\beta b_{2}a_{2}+2\beta b_{1}a_{1}-2\gamma  \lambda^2 b_{2}-2\gamma  \mu^2 a_{2}=0, \\
& \Big(\frac{G'}{G^2}\Big)^{1}&:& \ -k^2a_{1}+2\beta b_{1}a_{2}+2\beta a_{1}a_{0}-2\gamma \mu \lambda a_{1}=0, \\
& \Big(\frac{G'}{G^2}\Big)^{2}&:& \ -k^2a_{2}+\beta a_{1}^2+2 \beta a_{2}a_{0}-8\gamma \mu \lambda a_{2}=0, \\
& \Big(\frac{G'}{G^2}\Big)^{3}&:& \ \ \ \ 2 \beta a_{1}a_{2}-2 \gamma \lambda^2a_{1}=0, \\
& \Big(\frac{G'}{G^2}\Big)^{4} &:& \ \ \ \ \beta a_{2}^2-6\gamma \lambda^2a_{2}=0.
\end{aligned}
\end{equation*}

\noindent
Solving this system of equations through Maple, we get the following results;
\vspace{0.3cm}

\noindent
\textbf{Set1}: $k=\mp 2\sqrt{\lambda\gamma\mu}, \ \ \ \ a_{0}= \frac{6\lambda\gamma\mu}{\beta}, \ \ a_{1}=0, \ \ a_{2}= \frac{6\lambda^2\gamma}{\beta},\ \ b_{1}=0,\ \ b_{2}=0$,

\noindent
\textbf{Set2}: $k=\mp 2\sqrt{\lambda\gamma\mu}, \ \ \ \ a_{0}= \frac{6\lambda\gamma\mu}{\beta}, \ \ a_{1}=0, \ \ a_{2}=0,\ \ b_{1}=0,\ \ b_{2}=\frac{6\mu^2\gamma}{\beta}$,

\noindent
\textbf{Set3}: $k=\mp 2\sqrt{-\lambda\gamma\mu}, \ \ a_{0}= \frac{2\lambda\gamma\mu}{\beta}, \ \ a_{1}=0, \ \ a_{2}= \frac{6\lambda^2\gamma}{\beta},\ \ b_{1}=0,\ \ b_{2}=0$,

\noindent
\textbf{Set4}: $k=\mp 2\sqrt{-\lambda\gamma\mu}, \ \ a_{0}= \frac{2\lambda\gamma\mu}{\beta}, \ \ a_{1}=0, \ \ a_{2}=0,\ \ b_{1}=0,\ \ b_{2}=\frac{6\mu^2\gamma}{\beta}$,

\noindent
\textbf{Set5}: $k=\mp 4\sqrt{\lambda\gamma\mu}, \ \ \ \ a_{0}= \frac{12\lambda\gamma\mu}{\beta}, \ \ a_{1}=0, \ \ a_{2}= \frac{6\lambda^2\gamma}{\beta},\ \ b_{1}=0,\ \ b_{2}=\frac{6\mu^2\gamma}{\beta}$,

\noindent
\textbf{Set6}: $k=\mp 4\sqrt{-\lambda\gamma\mu}, \ \ a_{0}= \frac{-4\lambda\gamma\mu}{\beta}, \ \ a_{1}=0, \ \ a_{2}= \frac{6\lambda^2\gamma}{\beta},\ \ b_{1}=0,\ \ b_{2}=\frac{6\mu^2\gamma}{\beta}$.

\vspace{0.3cm}
\noindent
The above set of values gives the following exact solutions for conformable time fractional coupled Boussinesq equation.

\vspace{0.3cm}
\textbf{Solution 1}: 

\vspace{0.2cm}
\noindent
(i) If $\lambda\mu>0$, the trigonometric solution is found 
\begin{equation*}
\begin{aligned}
u_{1}(x,t)=U_{11}(\varepsilon)=\frac{6\lambda\gamma\mu}{\beta}+\frac{6\lambda^2\gamma}{\beta}\Bigg(\sqrt{\frac{\mu}{\lambda}}\Bigg(\frac{C\cos(\sqrt{\mu\lambda}\varepsilon)+D\sin(\sqrt{\mu\lambda}\varepsilon)}{D\cos(\sqrt{\mu\lambda}\varepsilon)-C\sin(\sqrt{\mu\lambda}\varepsilon)}\Bigg)\Bigg)^2.
\end{aligned}
\end{equation*}
From Eq.\eqref{eq14},
\begin{equation*}
\begin{aligned}
v_{1}(x,t)=-k \Bigg[\frac{6\lambda\gamma\mu}{\beta}+\frac{6\lambda^2\gamma}{\beta}\Bigg(\sqrt{\frac{\mu}{\lambda}}\Bigg(\frac{C\cos(\sqrt{\mu\lambda}\varepsilon)+D\sin(\sqrt{\mu\lambda}\varepsilon)}{D\cos(\sqrt{\mu\lambda}\varepsilon)-C\sin(\sqrt{\mu\lambda}\varepsilon)}\Bigg)\Bigg)^2\Bigg].
\end{aligned}
\end{equation*}

\noindent
(ii) If $\lambda\mu<0$, the hyperbolic solution is obtained
\begin{equation*}
\begin{aligned}
u_{2}(x,t)=U_{12}(\varepsilon)=\frac{6\lambda\gamma\mu}{\beta}+\frac{6\lambda^2\gamma}{\beta}\Bigg(-\frac{\sqrt{|\mu\lambda|}}{\lambda}\Bigg(\frac{C\sinh(2\sqrt{|\mu\lambda|}\varepsilon)+C\cosh(2\sqrt{|\mu\lambda|}\varepsilon)+D}{C\sinh(2\sqrt{|\mu\lambda|}\varepsilon)+C\cosh(2\sqrt{|\mu\lambda|}\varepsilon)-D}\Bigg)\Bigg)^2.
\end{aligned}
\end{equation*}
From Eq.\eqref{eq14},
\begin{equation*}
\begin{aligned}
v_{2}(x,t)=-k \Bigg[\frac{6\lambda\gamma\mu}{\beta}+\frac{6\lambda^2\gamma}{\beta}\Bigg(-\frac{\sqrt{|\mu\lambda|}}{\lambda}\Bigg(\frac{C\sinh(2\sqrt{|\mu\lambda|}\varepsilon)+C\cosh(2\sqrt{|\mu\lambda|}\varepsilon)+D}{C\sinh(2\sqrt{|\mu\lambda|}\varepsilon)+C\cosh(2\sqrt{|\mu\lambda|}\varepsilon)-D}\Bigg)\Bigg)^2\Bigg].
\end{aligned}
\end{equation*}

\noindent
(iii) If $\lambda\neq0,\ \ \ \mu=0$, the rational solution is found
\begin{equation*}
\begin{aligned}
u_{3}(x,t)=U_{13}(\varepsilon)=\frac{6\gamma}{\beta}\Bigg(\frac{C}{C\varepsilon+D}\Bigg)^2.
\end{aligned}
\end{equation*}
From Eq.\eqref{eq14},
\begin{equation*}
\begin{aligned}
v_{3}(x,t)=-k \Bigg[\frac{6\gamma}{\beta}\Bigg(\frac{C}{C\varepsilon+D}\Bigg)^2\Bigg],
\end{aligned}
\end{equation*}

\noindent
where $\varepsilon=x \mp 2\sqrt{\lambda\gamma\mu}\big(\frac{t^{\alpha}}{\alpha}\big) $.

\vspace{0.3cm}
\textbf{Solution 2}: 

\vspace{0.2cm}
\noindent
(i) If $\lambda\mu>0$, the trigonometric solution is found 
\begin{equation*}
\begin{aligned}
u_{4}(x,t)=U_{21}(\varepsilon)=\frac{6\lambda\gamma\mu}{\beta}+\frac{6\mu^2\gamma}{\beta}\Bigg(\sqrt{\frac{\mu}{\lambda}}\Bigg(\frac{C\cos(\sqrt{\mu\lambda}\varepsilon)+D\sin(\sqrt{\mu\lambda}\varepsilon)}{D\cos(\sqrt{\mu\lambda}\varepsilon)-C\sin(\sqrt{\mu\lambda}\varepsilon)}\Bigg)\Bigg)^{-2}.
\end{aligned}
\end{equation*}
From Eq.\eqref{eq14},
\begin{equation*}
\begin{aligned}
v_{4}(x,t)=-k \Bigg[\frac{6\lambda\gamma\mu}{\beta}+\frac{6\mu^2\gamma}{\beta}\Bigg(\sqrt{\frac{\mu}{\lambda}}\Bigg(\frac{C\cos(\sqrt{\mu\lambda}\varepsilon)+D\sin(\sqrt{\mu\lambda}\varepsilon)}{D\cos(\sqrt{\mu\lambda}\varepsilon)-C\sin(\sqrt{\mu\lambda}\varepsilon)}\Bigg)\Bigg)^{-2}\Bigg].
\end{aligned}
\end{equation*}
\noindent
(ii) If $\lambda\mu<0$, the hyperbolic solution is obtained
\begin{equation*}
\begin{aligned}
u_{5}(x,t)=U_{22}(\varepsilon)=\frac{6\lambda\gamma\mu}{\beta}+\frac{6\mu^2\gamma}{\beta}\Bigg(-\frac{\sqrt{|\mu\lambda|}}{\lambda}\Bigg(\frac{C\sinh(2\sqrt{|\mu\lambda|}\varepsilon)+C\cosh(2\sqrt{|\mu\lambda|}\varepsilon)+D}{C\sinh(2\sqrt{|\mu\lambda|}\varepsilon)+C\cosh(2\sqrt{|\mu\lambda|}\varepsilon)-D}\Bigg)\Bigg)^{-2}.
\end{aligned}
\end{equation*}
From Eq.\eqref{eq14},
\begin{equation*}
\begin{aligned}
v_{5}(x,t)=-k\Bigg[\frac{6\lambda\gamma\mu}{\beta}+\frac{6\mu^2\gamma}{\beta}\Bigg(-\frac{\sqrt{|\mu\lambda|}}{\lambda}\Bigg(\frac{C\sinh(2\sqrt{|\mu\lambda|}\varepsilon)+C\cosh(2\sqrt{|\mu\lambda|}\varepsilon)+D}{C\sinh(2\sqrt{|\mu\lambda|}\varepsilon)+C\cosh(2\sqrt{|\mu\lambda|}\varepsilon)-D}\Bigg)\Bigg)^{-2} \Bigg],
\end{aligned}
\end{equation*}
\noindent
where $\varepsilon=x \mp 2\sqrt{\lambda\gamma\mu}\big(\frac{t^{\alpha}}{\alpha}\big) $.

\vspace{0.3cm}
\textbf{Solution 3}: 

\vspace{0.2cm}
\noindent
(i) If $\lambda\mu>0$, the trigonometric solution is found 
\begin{equation*}
\begin{aligned}
u_{6}(x,t)=U_{31}(\varepsilon)=\frac{2\lambda\gamma\mu}{\beta}+\frac{6\lambda^2\gamma}{\beta}\Bigg(\sqrt{\frac{\mu}{\lambda}}\Bigg(\frac{C\cos(\sqrt{\mu\lambda}\varepsilon)+D\sin(\sqrt{\mu\lambda}\varepsilon)}{D\cos(\sqrt{\mu\lambda}\varepsilon)-C\sin(\sqrt{\mu\lambda}\varepsilon)}\Bigg)\Bigg)^2.
\end{aligned}
\end{equation*}
From Eq.\eqref{eq14},

\begin{equation*}
\begin{aligned}
v_{6}(x,t)=-k \Bigg[\frac{2\lambda\gamma\mu}{\beta}+\frac{6\lambda^2\gamma}{\beta}\Bigg(\sqrt{\frac{\mu}{\lambda}}\Bigg(\frac{C\cos(\sqrt{\mu\lambda}\varepsilon)+D\sin(\sqrt{\mu\lambda}\varepsilon)}{D\cos(\sqrt{\mu\lambda}\varepsilon)-C\sin(\sqrt{\mu\lambda}\varepsilon)}\Bigg)\Bigg)^2\Bigg].
\end{aligned}
\end{equation*}
\noindent
(ii) If $\lambda\mu<0$, the hyperbolic solution is obtained
\begin{equation*}
\begin{aligned}
& u_{7}(x,t)=U_{32}(\varepsilon)=\frac{2\lambda\gamma\mu}{\beta}+\frac{6\lambda^2\gamma}{\beta}\Bigg(-\frac{\sqrt{|\mu\lambda|}}{\lambda}\Bigg(\frac{C\sinh(2\sqrt{|\mu\lambda|}\varepsilon)+C\cosh(2\sqrt{|\mu\lambda|}\varepsilon)+D}{C\sinh(2\sqrt{|\mu\lambda|}\varepsilon)+C\cosh(2\sqrt{|\mu\lambda|}\varepsilon)-D}\Bigg)\Bigg)^2.
\end{aligned}
\end{equation*}
From Eq.\eqref{eq14},

\begin{equation*}
\begin{aligned}
v_{7}(x,t)=-k \Bigg[\frac{2\lambda\gamma\mu}{\beta}+\frac{6\lambda^2\gamma}{\beta}\Bigg(-\frac{\sqrt{|\mu\lambda|}}{\lambda}\Bigg(\frac{C\sinh(2\sqrt{|\mu\lambda|}\varepsilon)+C\cosh(2\sqrt{|\mu\lambda|}\varepsilon)+D}{C\sinh(2\sqrt{|\mu\lambda|}\varepsilon)+C\cosh(2\sqrt{|\mu\lambda|}\varepsilon)-D}\Bigg)\Bigg)^2 \Bigg].
\end{aligned}
\end{equation*}

\noindent
(iii) If $\lambda\neq0,\ \ \ \mu=0$, the rational solution is found
\begin{equation*}
\begin{aligned}
u_{8}(x,t)=U_{33}(\varepsilon)=\frac{6\gamma}{\beta}\Bigg(\frac{C}{C\varepsilon+D}\Bigg)^2.
\end{aligned}
\end{equation*}
From Eq.\eqref{eq14},
\begin{equation*}
\begin{aligned}
v_{8}(x,t)=-k \Bigg[\frac{6\gamma}{\beta}\Bigg(\frac{C}{C\varepsilon+D}\Bigg)^2 \Bigg],
\end{aligned}
\end{equation*}
\noindent
where $\varepsilon=x \mp 2\sqrt{-\lambda\gamma\mu}\big(\frac{t^{\alpha}}{\alpha}\big) $.

\vspace{0.3cm}
\textbf{Solution 4}:

\vspace{0.2cm}
\noindent
(i) If $\lambda\mu>0$, the trigonometric solution is found 
\begin{equation*}
\begin{aligned}
u_{9}(x,t)=U_{41}(\varepsilon)=\frac{2\lambda\gamma\mu}{\beta}+\frac{6\mu^2\gamma}{\beta}\Bigg(\sqrt{\frac{\mu}{\lambda}}\Bigg(\frac{C\cos(\sqrt{\mu\lambda}\varepsilon)+D\sin(\sqrt{\mu\lambda}\varepsilon)}{D\cos(\sqrt{\mu\lambda}\varepsilon)-C\sin(\sqrt{\mu\lambda}\varepsilon)}\Bigg)\Bigg)^{-2}.
\end{aligned}
\end{equation*}
From Eq.\eqref{eq14},

\begin{equation*}
\begin{aligned}
v_{9}(x,t)=-k \Bigg[\frac{2\lambda\gamma\mu}{\beta}+\frac{6\mu^2\gamma}{\beta}\Bigg(\sqrt{\frac{\mu}{\lambda}}\Bigg(\frac{C\cos(\sqrt{\mu\lambda}\varepsilon)+D\sin(\sqrt{\mu\lambda}\varepsilon)}{D\cos(\sqrt{\mu\lambda}\varepsilon)-C\sin(\sqrt{\mu\lambda}\varepsilon)}\Bigg)\Bigg)^{-2} \Bigg].
\end{aligned}
\end{equation*}

\noindent
(ii) If $\lambda\mu<0$, the hyperbolic solution is obtained
\begin{equation*}
\begin{aligned}
u_{10}(x,t)=U_{42}(\varepsilon)=\frac{2\lambda\gamma\mu}{\beta}+\frac{6\mu^2\gamma}{\beta}\Bigg(-\frac{\sqrt{|\mu\lambda|}}{\lambda}\Bigg(\frac{C\sinh(2\sqrt{|\mu\lambda|}\varepsilon)+C\cosh(2\sqrt{|\mu\lambda|}\varepsilon)+D}{C\sinh(2\sqrt{|\mu\lambda|}\varepsilon)+C\cosh(2\sqrt{|\mu\lambda|}\varepsilon)-D}\Bigg)\Bigg)^{-2}.
\end{aligned}
\end{equation*}
From Eq.\eqref{eq14},
\begin{equation*}
\begin{aligned}
v_{10}(x,t)=-k \Bigg[ \frac{2\lambda\gamma\mu}{\beta}+\frac{6\mu^2\gamma}{\beta}\Bigg(-\frac{\sqrt{|\mu\lambda|}}{\lambda}\Bigg(\frac{C\sinh(2\sqrt{|\mu\lambda|}\varepsilon)+C\cosh(2\sqrt{|\mu\lambda|}\varepsilon)+D}{C\sinh(2\sqrt{|\mu\lambda|}\varepsilon)+C\cosh(2\sqrt{|\mu\lambda|}\varepsilon)-D}\Bigg)\Bigg)^{-2}\Bigg],
\end{aligned}
\end{equation*}
\noindent
where $\varepsilon=x \mp 2\sqrt{-\lambda\gamma\mu}\big(\frac{t^{\alpha}}{\alpha}\big) $.

\vspace{0.3cm}
\textbf{Solution 5}: 

\vspace{0.2cm}
\noindent
(i) If $\lambda\mu>0$, the trigonometric solution is found 
\begin{equation*}
\begin{aligned}
& u_{11}(x,t)=U_{51}(\varepsilon)=\frac{12\lambda\gamma\mu}{\beta}+\frac{6\lambda^2\gamma}{\beta}\Bigg(\sqrt{\frac{\mu}{\lambda}}\Bigg(\frac{C\cos(\sqrt{\mu\lambda}\varepsilon)+D\sin(\sqrt{\mu\lambda}\varepsilon)}{D\cos(\sqrt{\mu\lambda}\varepsilon)-C\sin(\sqrt{\mu\lambda}\varepsilon)}\Bigg)\Bigg)^2 \\
& \ \ \ \ \ \ \ \ \ \ \ \ \ \ \ \ \ \ \ \ \ \ \ \ \ \ \ \ \ \ \ \ \ +\frac{6\mu^2\gamma}{\beta}\Bigg(\sqrt{\frac{\mu}{\lambda}}\Bigg(\frac{C\cos(\sqrt{\mu\lambda}\varepsilon)+D\sin(\sqrt{\mu\lambda}\varepsilon)}{D\cos(\sqrt{\mu\lambda}\varepsilon)-C\sin(\sqrt{\mu\lambda}\varepsilon)}\Bigg)\Bigg)^{-2}.
\end{aligned}
\end{equation*}
From Eq.\eqref{eq14},
\begin{equation*}
\begin{aligned}
&v_{11}(x,t)=-k \Bigg[ \frac{12\lambda\gamma\mu}{\beta}+\frac{6\lambda^2\gamma}{\beta}\Bigg(\sqrt{\frac{\mu}{\lambda}}\Bigg(\frac{C\cos(\sqrt{\mu\lambda}\varepsilon)+D\sin(\sqrt{\mu\lambda}\varepsilon)}{D\cos(\sqrt{\mu\lambda}\varepsilon)-C\sin(\sqrt{\mu\lambda}\varepsilon)}\Bigg)\Bigg)^2 \\
& \ \ \ \ \ \ \ \ \ \ \ \ \ \ \ \ \ \ \ \ \ \ \ \ \ \ \ \ \ \ \ \ \ +\frac{6\mu^2\gamma}{\beta}\Bigg(\sqrt{\frac{\mu}{\lambda}}\Bigg(\frac{C\cos(\sqrt{\mu\lambda}\varepsilon)+D\sin(\sqrt{\mu\lambda}\varepsilon)}{D\cos(\sqrt{\mu\lambda}\varepsilon)-C\sin(\sqrt{\mu\lambda}\varepsilon)}\Bigg)\Bigg)^{-2}\Bigg].
\end{aligned}
\end{equation*}
\noindent
(ii) If $\lambda\mu<0$, the hyperbolic solution is obtained
\begin{equation*}
\begin{aligned}
& u_{12}(x,t)=U_{52}(\varepsilon)=\frac{12\lambda\gamma\mu}{\beta}+\frac{6\lambda^2\gamma}{\beta}\Bigg(-\frac{\sqrt{|\mu\lambda|}}{\lambda}\Bigg(\frac{C\sinh(2\sqrt{|\mu\lambda|}\varepsilon)+C\cosh(2\sqrt{|\mu\lambda|}\varepsilon)+D}{C\sinh(2\sqrt{|\mu\lambda|}\varepsilon)+C\cosh(2\sqrt{|\mu\lambda|}\varepsilon)-D}\Bigg)\Bigg)^2 \\
& \ \ \ \ \ \ \ \ \ \ \ \ \ \ \ \ \ \ \ \ \ \ \ \ \ \ \ \ \ \ \ \ \ +\frac{6\mu^2\gamma}{\beta}\Bigg(-\frac{\sqrt{|\mu\lambda|}}{\lambda}\Bigg(-\frac{C\sinh(2\sqrt{|\mu\lambda|}\varepsilon)+C\cosh(2\sqrt{|\mu\lambda|}\varepsilon)+D}{C\sinh(2\sqrt{|\mu\lambda|}\varepsilon)+C\cosh(2\sqrt{|\mu\lambda|}\varepsilon)-D}\Bigg)\Bigg)^{-2}. 
\end{aligned}
\end{equation*}
From Eq.\eqref{eq14},
\begin{equation*}
\begin{aligned}
& v_{12}(x,t)=-k \Bigg[\frac{12\lambda\gamma\mu}{\beta}+\frac{6\lambda^2\gamma}{\beta}\Bigg(-\frac{\sqrt{|\mu\lambda|}}{\lambda}\Bigg(\frac{C\sinh(2\sqrt{|\mu\lambda|}\varepsilon)+C\cosh(2\sqrt{|\mu\lambda|}\varepsilon)+D}{C\sinh(2\sqrt{|\mu\lambda|}\varepsilon)+C\cosh(2\sqrt{|\mu\lambda|}\varepsilon)-D}\Bigg)\Bigg)^2 \\
& \ \ \ \ \ \ \ \ \ \ \ \ \ \ \ \ \ \ \ \ \ \ \ \ \ \ \ \ \ \ \ \ \ +\frac{6\mu^2\gamma}{\beta}\Bigg(-\frac{\sqrt{|\mu\lambda|}}{\lambda}\Bigg(-\frac{C\sinh(2\sqrt{|\mu\lambda|}\varepsilon)+C\cosh(2\sqrt{|\mu\lambda|}\varepsilon)+D}{C\sinh(2\sqrt{|\mu\lambda|}\varepsilon)+C\cosh(2\sqrt{|\mu\lambda|}\varepsilon)-D}\Bigg)\Bigg)^{-2} \Bigg].
\end{aligned}
\end{equation*} 
\noindent
(iii) If $\lambda\neq0,\ \ \ \mu=0$, the rational solution is found
\begin{equation*}
\begin{aligned}
u_{13}(x,t)=U_{53}(\varepsilon)=\frac{6\gamma}{\beta}\Bigg(\frac{C}{C\varepsilon+D}\Bigg)^2.
\end{aligned}
\end{equation*}
From Eq.\eqref{eq14},
\begin{equation*}
\begin{aligned}
v_{13}(x,t)=-k \Bigg[\frac{6\gamma}{\beta}\Bigg(\frac{C}{C\varepsilon+D}\Bigg)^2\Bigg],
\end{aligned}
\end{equation*}
\noindent
where $\varepsilon=x \mp 4\sqrt{\lambda\gamma\mu}\big(\frac{t^{\alpha}}{\alpha}\big) $.

\vspace{0.3cm}
\textbf{Solution 6}: 

\vspace{0.2cm}
\noindent
(i) If $\lambda\mu>0$, the trigonometric solution is found 
\begin{equation*}
\begin{aligned}
& u_{14}(x,t)=U_{61}(\varepsilon)=\frac{-4\lambda\gamma\mu}{\beta}+\frac{6\lambda^2\gamma}{\beta}\Bigg(\sqrt{\frac{\mu}{\lambda}}\Bigg(\frac{C\cos(\sqrt{\mu\lambda}\varepsilon)+D\sin(\sqrt{\mu\lambda}\varepsilon)}{D\cos(\sqrt{\mu\lambda}\varepsilon)-C\sin(\sqrt{\mu\lambda}\varepsilon)}\Bigg)\Bigg)^2 \\
& \ \ \ \ \ \ \ \ \ \ \ \ \ \ \ \ \ \ \ \ \ \ \ \ \ \ \ \ \ \ \ \ \ +\frac{6\mu^2\gamma}{\beta}\Bigg(\sqrt{\frac{\mu}{\lambda}}\Bigg(\frac{C\cos(\sqrt{\mu\lambda}\varepsilon)+D\sin(\sqrt{\mu\lambda}\varepsilon)}{D\cos(\sqrt{\mu\lambda}\varepsilon)-C\sin(\sqrt{\mu\lambda}\varepsilon)}\Bigg)\Bigg)^{-2}.
\end{aligned}
\end{equation*}
From Eq.\eqref{eq14},
\begin{equation*}
\begin{aligned}
&v_{14}(x,t)=-k \Bigg[ \frac{-4\lambda\gamma\mu}{\beta}+\frac{6\lambda^2\gamma}{\beta}\Bigg(\sqrt{\frac{\mu}{\lambda}}\Bigg(\frac{C\cos(\sqrt{\mu\lambda}\varepsilon)+D\sin(\sqrt{\mu\lambda}\varepsilon)}{D\cos(\sqrt{\mu\lambda}\varepsilon)-C\sin(\sqrt{\mu\lambda}\varepsilon)}\Bigg)\Bigg)^2 \\
& \ \ \ \ \ \ \ \ \ \ \ \ \ \ \ \ \ \ \ \ \ \ \ \ \ \ \ \ \ \ \ \ \ +\frac{6\mu^2\gamma}{\beta}\Bigg(\sqrt{\frac{\mu}{\lambda}}\Bigg(\frac{C\cos(\sqrt{\mu\lambda}\varepsilon)+D\sin(\sqrt{\mu\lambda}\varepsilon)}{D\cos(\sqrt{\mu\lambda}\varepsilon)-C\sin(\sqrt{\mu\lambda}\varepsilon)}\Bigg)\Bigg)^{-2}\Bigg].
\end{aligned}
\end{equation*}

\noindent
(ii) If $\lambda\mu<0$, the hyperbolic solution is obtained
\begin{equation*}
\begin{aligned}
& u_{15}(x,t)=U_{62}(\varepsilon)=\frac{-4\lambda\gamma\mu}{\beta}+\frac{6\lambda^2\gamma}{\beta}\Bigg(-\frac{\sqrt{|\mu\lambda|}}{\lambda}\Bigg(\frac{C\sinh(2\sqrt{|\mu\lambda|}\varepsilon)+C\cosh(2\sqrt{|\mu\lambda|}\varepsilon)+D}{C\sinh(2\sqrt{|\mu\lambda|}\varepsilon)+C\cosh(2\sqrt{|\mu\lambda|}\varepsilon)-D}\Bigg)\Bigg)^2 \\
& \ \ \ \ \ \ \ \ \ \ \ \ \ \ \ \ \ \ \ \ \ \ \ \ \ \ \ \ \ \ \ \ \ +\frac{6\mu^2\gamma}{\beta}\Bigg(-\frac{\sqrt{|\mu\lambda|}}{\lambda}\Bigg(-\frac{C\sinh(2\sqrt{|\mu\lambda|}\varepsilon)+C\cosh(2\sqrt{|\mu\lambda|}\varepsilon)+D}{C\sinh(2\sqrt{|\mu\lambda|}\varepsilon)+C\cosh(2\sqrt{|\mu\lambda|}\varepsilon)-D}\Bigg)\Bigg)^{-2}.
\end{aligned}
\end{equation*}
From Eq.\eqref{eq14},
\begin{equation*}
\begin{aligned}
&v_{15}(x,t)=-k \Bigg[ \frac{-4\lambda\gamma\mu}{\beta}+\frac{6\lambda^2\gamma}{\beta}\Bigg(-\frac{\sqrt{|\mu\lambda|}}{\lambda}\Bigg(\frac{C\sinh(2\sqrt{|\mu\lambda|}\varepsilon)+C\cosh(2\sqrt{|\mu\lambda|}\varepsilon)+D}{C\sinh(2\sqrt{|\mu\lambda|}\varepsilon)+C\cosh(2\sqrt{|\mu\lambda|}\varepsilon)-D}\Bigg)\Bigg)^2 \\
& \ \ \ \ \ \ \ \ \ \ \ \ \ \ \ \ \ \ \ \ \ \ \ \ \ \ \ \ \ \ \ \ \ +\frac{6\mu^2\gamma}{\beta}\Bigg(-\frac{\sqrt{|\mu\lambda|}}{\lambda}\Bigg(-\frac{C\sinh(2\sqrt{|\mu\lambda|}\varepsilon)+C\cosh(2\sqrt{|\mu\lambda|}\varepsilon)+D}{C\sinh(2\sqrt{|\mu\lambda|}\varepsilon)+C\cosh(2\sqrt{|\mu\lambda|}\varepsilon)-D}\Bigg)\Bigg)^{-2}\Bigg].
\end{aligned}
\end{equation*}

\noindent
(iii) If $\lambda\neq0,\ \ \ \mu=0$, the rational solution is found
\begin{equation*}
\begin{aligned}
u_{16}(x,t)=U_{63}(\varepsilon)=\frac{6\gamma}{\beta}\Bigg(\frac{C}{C\varepsilon+D}\Bigg)^2.
\end{aligned}
\end{equation*}
From Eq.\eqref{eq14},
\begin{equation*}
\begin{aligned}
 v_{16}(x,t)=-k \Bigg[\frac{6\gamma}{\beta}\Bigg(\frac{C}{C\varepsilon+D}\Bigg)^2\Bigg],
\end{aligned}
\end{equation*}

\noindent
where $\varepsilon=x  \mp 4\sqrt{-\lambda\gamma\mu}\big(\frac{t^{\alpha}}{\alpha}\big) $.

\vspace{0.2cm}

\noindent
In Figure 2, the physical characteristics of Eq.(4.5) whose solutions are used in \textbf{Solution 3} have been shown for $\alpha = 0.5$ with some special values.
\begin{figure}
\centering
  \includegraphics[width=14cm]{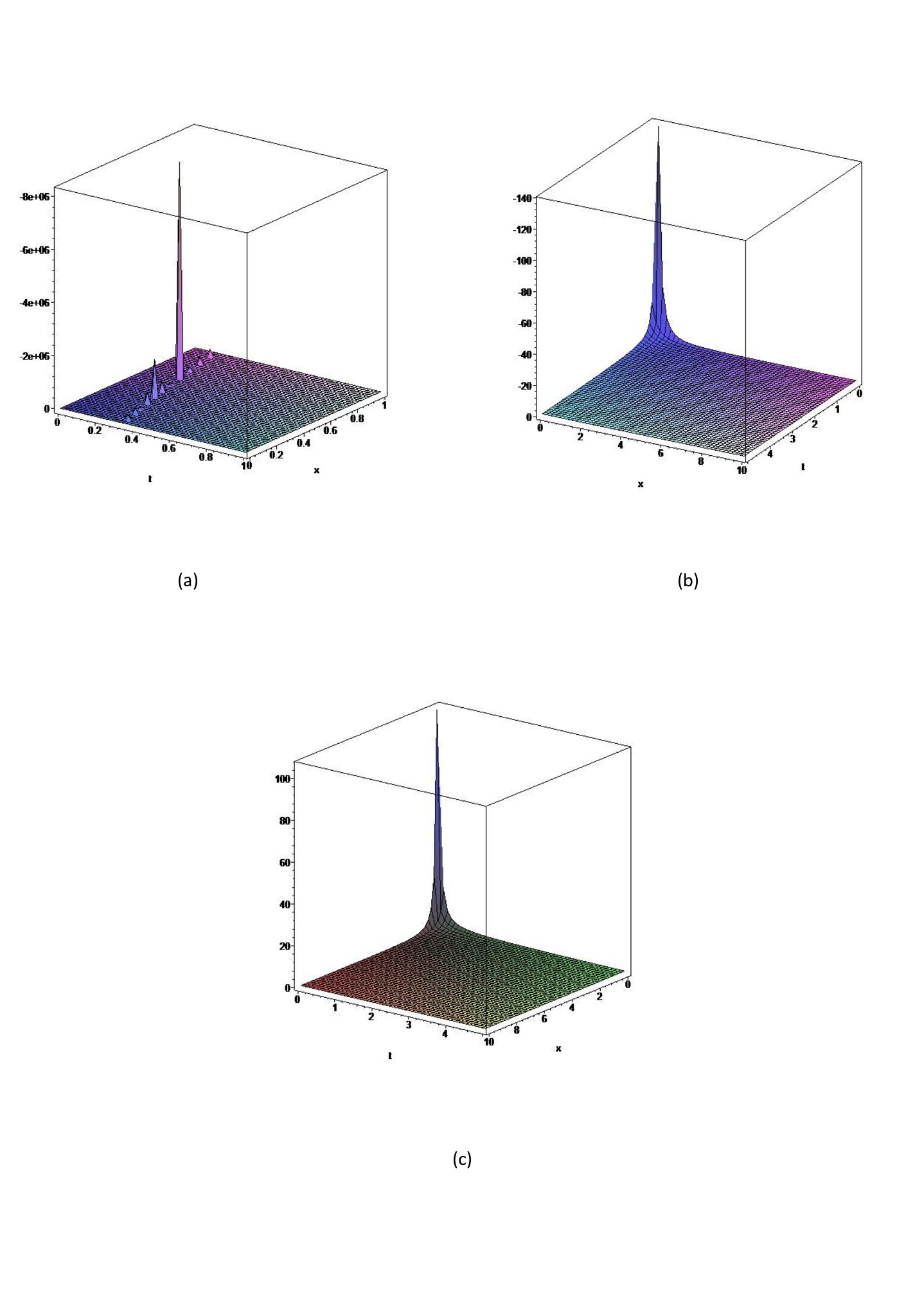}\\
  \caption{(a) The trigonometric solution of $u(x,t)$ when $ \lambda=\mu=\beta=C=D=1, \gamma=-1.$ (b) The hyperbolic solution of $u(x,t)$ when $ \lambda=0.5, \mu=-0.3, C=D=\beta=\gamma=1$. (c) The rational solution of $u(x,t)$ when $ \lambda=1, \mu=0, C=D=\beta=1,\gamma=-1$.}
\end{figure}

\section{Conlusion}
In this paper, analytical solutions of conformable time-fractional nonlinear Boussinesq equations which known as an important class of fractional differential equations in mathematical physics were found by $\Big(\frac{G'}{G^2}\Big)$-expansion method. Different types of solutions were officially found, including trigonometric, hyperbolic and rational function solutions. The obtained solutions were not the same and  they were much recent than previously results. It was ensured that the solutions found were correct by replacing the original equations. Since $\Big(\frac{G'}{G^2}\Big)$-expansion method is efficient and not complicated, it can also be used to get solutions for other fractional differential equations encountered in mathematical physics. 

%% main text
%\input{chapters/en/Introduction}
%
%\input{chapters/en/Methodology}
%
%\input{chapters/en/Data}
%
%\input{chapters/en/Results}
%
%\input{chapters/en/Conclusion}
%
%\input{chapters/en/Acknowledgement}
%
%\appendix
%\input{chapters/en/Appendix}

\end{document}